\documentclass[10pt]{amsart}
\usepackage{amsmath,amssymb,amsfonts}
\usepackage{enumerate}
\usepackage{amscd, amssymb, latexsym, amsmath, amscd}
\usepackage[all]{xy}

\newtheorem{theorem}{Theorem}[section]
\newtheorem{thm}[theorem]{Theorem}
\newtheorem{prop}[theorem]{Proposition}
\newtheorem{cor}[theorem]{Corollary}
\newtheorem{conj}[theorem]{Conjecture}
\newtheorem{lemma}[theorem]{Lemma}
\newtheorem{question}[theorem]{Question}

\newcommand{\R}{\mathbb R}
\newcommand{\C}{\mathbb C}
\newcommand{\B}{\mathbb B}
\newcommand{\D}{\mathbb D}

\newcommand{\A}{\mathbb A}

\def\SL{{\rm SL}}

\def\GL{{\rm GL}}

\def\O{{\rm O}}

\makeatletter
\def\l@section{\@tocline{1}{4pt}{1pc}{}{}}
\def\l@subsection{\@tocline{2}{0pt}{2pc}{5pc}{}}
\makeatother
\begin{document}

\title{Self-dual representations of division algebras and Weil groups: A contrast}
\author{Dipendra Prasad  and Dinakar Ramakrishnan}
\thanks{D.~Prasad was partially supported by Friends of the IAS and the Von Neumann Fund, and then by the Clay Math Institute, while D.~Ramakrishnan
was partially supported by the National Science Foundation}
\date{}
\maketitle
\pagestyle{myheadings}
\markboth{D.~Prasad and D.~Ramakrishnan}{Self-dual representations}

\section*{\bf Introduction}

If $\rho$ is a selfdual representation of a group $G$ on a vector
space $V$ over ${\Bbb C}$,
we will say that $\rho$ is {\it orthogonal}, resp. {\it symplectic},
if $G$ leaves a nondegenerate
{\it symmetric}, resp. {\it alternating}, bilinear
form $B$ on $V$ invariant. If $\rho$ is irreducible, exactly one of
these possibilities will occur, with $B$ unique up to scaling,
and we may define a {\it sign} $c(\rho) \in \{\pm
1\}$, taken to be $+1$, resp. $-1$, in the orthogonal, resp.
symplectic, case.

Now let $k$ be a local field of characteristic 0.
The groups of interest to us will be $G=$\GL$_m(D)$,
where $D$ is a division algebra of center
$k$ and index $r$, and the Weil group $W_k$. Note that the group GL$_1(D)$ ($=D^\times$) is
even compact modulo the center, hence its complex irreducible
representations $\pi$ are finite dimensional.
For any $m$, let
$R^0(G)$ denote the set of irreducible admissible $\C$-representations $\pi$ (up to
equivalence) which correspond, by the generalized Jacquet-Langlands
correspondence (\cite{bad}), to irreducible discrete series
representations $\pi^\prime$ of GL$_n(k)$, with $n=mr$.
When used in conjunction with the local Langlands correspondence
(\cite{HaT}, \cite{He}), there exists a bijection  $\pi\mapsto\sigma$ satisfying certain
natural properties, such as the preservation of
$\varepsilon$-factors of pairs, between $R^0(G)$ and the set $Irr_n(W^\prime_k)$
of irreducible representations $\sigma$ of
$W_k^\prime$ of dimension $n$, again taken up to equivalence. Here $W_k^\prime$
denotes $W_k$ if $k$ is archimedean, and the extended Weil group
$W_k \times {\rm SL}(2,\C)$ if $k$ is non-archimedean.
One calls $\sigma$ the {\it Langlands
parameter} of $\pi$. It is immediate from the construction that
$\pi$ is selfdual if and only if  $\sigma$ is. However, the local
Langlands reciprocity is not {\it a priori} sensitive to the {\it
finer question} of whether $c(\pi)$ equals $c(\sigma)$ or
$-c(\sigma)$. The main result of this paper is
the following.

\medskip

\noindent{\bf Theorem A} \, \it Let $n=mr$, $D$ a division algebra of
index $r$ over a local field $k$ of characteristic zero, $G={\rm GL}_m(D)$,
and $\pi$ a selfdual representation of $G$ in $R^0(G)$ with parameter
$\sigma \in Irr_n(W^\prime_k)$. Then we have
$$
(-1)^mc(\pi) \, = \, (-1)^nc(\sigma)^m.
$$
\rm

\medskip

\noindent{\bf Corollary B} \, \it Let $\pi$ be an irreducible, selfdual representation of $D^\times$,
for any division algebra of index $n$
over a local field $k$ of characteristic zero, and let $\sigma$ be the $n$-dimensional parameter of
$\pi$. If $n$ is odd, $\pi$ is
always orthogonal, while if $n$ is even,
$$
\pi \quad {\rm orthogonal} \quad \Longleftrightarrow \quad \sigma \quad {\rm symplectic}.
$$
\rm

\medskip

When $n$ is odd, $\sigma$ is necessarily orthogonal,
and so Theorem A implies, for any factorization $n=mr$ and $G=$ GL$_m(D)$ with $D$ a division
algebra of index $r$,
that $c(\pi)=+1$, i.e., $\pi$ is orthogonal like its parameter. Now let $m=1$.
Then we get $c(\pi)=(-1)^{n+1}c(\sigma)$, which
implies that {\it for $n$ even, $\pi$ is symplectic if and only if $\sigma$ is odd}, proving the Corollary, assuming the Theorem.
This surprising {\it flip} for $n$ even is what
we noticed first for $n=2$, spurring our interest in the general case, which is more subtle to establish.
Based on considerations of Poincar\'e
duality on the middle dimensional cohomology of certain coverings of
the Drinfeld upper-half space, we conjectured in \cite{PR} the assertion of Corollary B,
and established  some positive results in \cite{PR} and \cite{P1}, including
the case of $n=2$. In \cite{P1} it was proved that if
$n$ is odd and if the residual characteristic of $k$ is odd, then
$D^\times$ has no selfdual irreducible representations of  dimension $>
1$, showing that in this case, the conjecture is difficult only for the even
residual characteristic. An engaging program to prove Corollary B along the geometric lines,
using cohomological methods
involving the formal moduli of Lubin-Tate groups, has been announced in the supercuspidal case
by Laurent Fargues; it does not seem, however, that, without further input, his suggested methods would work for general discrete series representations,
nor for GL$_m(D)$ for $m>1$.
It may also be useful to note that, when $r=1$,
so that $G=$ GL$_n(k)$, the formula of Theorem A is easily seen to yield $c(\pi)=+1$, which is elementary to
prove directly,
but is nevertheless useful while using global arguments to deal with more
 difficult cases.

Let $n$ be even and $m=1$. Then Corollary B associates, to each irreducible, symplectic
Galois representation $\sigma$ of dimension $n$, a
{\it new secondary invariant},
defined by whether, or not, the associated orthogonal representation $\pi$ of $D^\times$ lifts
to the (s)pin group. This aspect
was investigated for $n=2$ in \cite{PR}.

\medskip

Our proof of Theorem A for non-archimedean $k$ proceeds by using {\it global methods}, and the philosophy is that if
the sign $c(\pi)$ is the asserted one for local representations $\pi$ of simple type at one prime, then the sign
should also be right for any complicated local representation at a {\it different} prime. This is why
section 4 of the paper is spent verifying the assertion in a special case, namely when
$G=D^\times$, $n$ is even, and $\pi$ is of {\sl level
$1$}, i.e., trivial on the level $1$ subgroup $D^\times(1)$. When the parameter is not trivial on SL$(2,\C)$,
we have to use a variant of the global approach, having to deal (in the process)
with generalized Steinberg representations of GL$(n,k)$.
When $\pi$ is not of this type, we globalize appropriately and get a level $1$ representation at another prime.
When the associated representation of GL$(n,k)$ is supercuspidal, We use the strong functoriality
transfer between generic automorphic representations of classical
groups and those of the corresponding general linear group admitting a pole of an
appropriate $L$-function, which have been established in a series of
papers due to many people, namely \cite{CoKPSS}, \cite{CoPSS},
\cite{J-S1}, \cite{GJR}, \cite{J-S2}). A synthesis of these results is a necessary first
step in proving Theorem A for $G=D^\times$, and even there the case when the parameter is non-trivial
on SL$(2,\C)$ requires a new twist, involving a different sort of globalization (see section 3), which
allows us to transport the sign problem to one with trivial SL$(2)$-parameter. We then use the result in this crucial case to
attack the general situation of $G=$ GL$_m(D)$ later.

\medskip

The workhorse which allows us to use global
methods is the following {\it product formula}, reminiscent of
quadratic reciprocity.

\medskip

\noindent{\bf Theorem C} \, \it Let $F$ be a global field,
$G=$GL$_m(\D)$, where $\D$ is a division algebra of dimension $r^2$
over $F$ and $Z$ the center of $G$. Suppose $\Pi = \otimes_v^\prime
\Pi_v$ is an irreducible, selfdual automorphic representation of
$G(\A_F)$ of central character $\omega$, which occurs with
multiplicity one in $L^2(G(F)Z(\A_F)\backslash G(\A_F),\omega)$.
Then we have
$$
c(\Pi) \, = \, \prod\limits_{v \in {\rm ram}(\D)} \, c(\Pi_v) \, = \, 1,
$$
where ram$(\D)$ denotes the set of places where $\D$ is ramified.
\rm

\medskip

As a consequence, we see that in the case $m=1$ and
ram$(\D)=\{u,v\}$, we have $c(\Pi_u)=c(\Pi_v)$, i.e.,
$\Pi_u$ and $\Pi_v$ are both orthogonal or both symplectic. In
particular, if we know one, we know the other.

There are two reasons why Theorem C holds. The first is that the
{\it global sign} $c(\pi)$ is $1$, and this works in greater
generality by a simple argument (see section 2). The second is
that the local sign $c(\pi_v)$ is $1$ whenever $G(F_v)$ is
GL$_n(F_v)$.

\medskip

Thanks to the product formula, given an irreducible selfdual
representation $\pi$ of $D^\times$, with $D$ a division algebra over a
local field $k$, our strategy becomes one of trying to find a number
field $F$ with $F_v=k$ for a place $v$ of $F$, a division algebra
$\D$ over $F$ with ram$(\D)=\{u,v\}$, and a {\it selfdual}
automorphic representation $\Pi$ of $\D^\times(\A_F)$ such that
$\Pi_v=\pi$ and $\Pi_u$ a known representation. Clearly, the
parameters of $\pi$ and $\Pi_u$ must have the same parity. Thanks to
the generalized Jacquet-Langlands correspondence between $\D^\times$
and GL$_n/F$, we see that it suffices to find a {\it selfdual}
discrete automorphic representation $\Pi^\prime$ of GL$_n(\A_F)$
with $\Pi_v^\prime$, resp. $\Pi_u^\prime$, being associated to
$\Pi_v$, resp. $\Pi_u$. The difficulty is not in globalizing, but in
choosing a global $\Pi^\prime$ which is also selfdual. We have,
thankfully, been able to  address this, when $\Pi_u$ and $\Pi_v$ are
both supercuspidal by using various
known instances of Functoriality, culminating in the works of Jiang and Soudry
(\cite{J-S2}, \cite{Sou}).

\medskip

The statement of Theorem A for $n$ odd is the simplest
as it does not refer to the Langlands parameter at all. It says in particular that a
selfdual irreducible
representation of $D^\times$ for a division algebra of index $n$, an odd integer, must have a
non-degenerate invariant symmetric bilinear form. Since $D^\times/k^\times$ is a
compact, profinite group, the question is clearly in the realm of finite group
theory. However, our proof  uses many recent and nontrivial
results in the theory of
Automorphic representations to achieve this. Recently, Bushnell
and Henniart have given a local proof of this case of our result in \cite{BH}.

\medskip

More generally, given a connected reductive algebraic group $G/k$ with an
involution $\theta$, and an irreducible $\theta$-selfdual
representation $\pi$ of $G(k)$, we associate, in section 7, an
invariant $c_\theta(\pi)\in\{\pm 1\}$. When $\theta = 1$, $c_\theta(\pi)$ gives
exactly the information on whether $\pi$ is orthogonal or
symplectic. There is an analogous invariant for Galois
representations in chapter 15 of \cite{Rog}.
We would like to know if $c_\theta(\pi)$ is the same for two such
representations in the same $L$-packet, and if so, whether $c_\theta(\pi)$
admits a formula in terms of the Langlands parameter. In any case,
our proof of Theorem C in section 2 extends to this situation when one has
multiplicity one for the relevant global representations. We also show locally that for
suitable $\theta$, generic discrete series representations $\pi$ of
a quasi-split $G$ have $c_\theta(\pi)=1$.

When the dual group ${}^LG$ admits an involution ${}^L\theta$, we define
(in section 8) an invariant $c_{_{{}^L\theta}}(\sigma)\in \{\pm 1\}$
for $L$-homomorphisms $\sigma$ from $W_k^\prime$ into ${}^LG$ which are equivalent
to $\sigma^\theta$ under the connected component of ${}^LG$.
When $G/k$ is a quasi-split reductive group with an involution $\theta$
stabilizing a Borel subgroup, its maximal torus, and the
associated pinned root system, one gets such a
dual involution ${}^L\theta$ of ${}^LG$. It will be interesting to make a comparison
of $c_{_{}^L\theta}(\sigma)$, for $\sigma$ whose image in ${}^LG$
is not contained in a Levi subgroup,
and the $c_\theta(\pi)$ for $\pi$, for $\pi$ in the $L$-packet associated to $\sigma$ by the (conjectural) local
Langlands correspondence.

\medskip

The first author would like to thank the California Institute of
Technology for the invitation to visit, enabling the authors to
collaborate on this work. He thanks the Institute for Advanced Study where
a significant portion of
the work was done, and gratefully acknowledges receiving support through
 grants to the Institute by the Friends of
the Institute, and the von Neumann Fund, and also thanks the Clay Math Institute for supporting him at
the final stages of this work. He would also like to thank
E.~Lapid for his interest and for suggesting a similar approach for
$G=D^\times$. Both the authors would like to thank J.~Bernstein, for inquiring about GL$_m(D)$ for $m>1$,
as well as L.~Clozel, M.~Harris,
G.~Henniart, and W.~ Zink for making helpful comments. They also
want to express their thanks to Dihua Jiang and David Soudry for sending them a letter
(\cite{Sou-lett}) explaining the completion of the proof of their theorem concerning the strong
generic transfer between classical groups and GL$(n)$; we summarize their
results in section 2. The second author would in addition
like to thank Herv\'e Jacquet and Akshay Venkatesh for their
interest, and also acknowledge support during the course of this work from the National Science
Foundation through the grants DMS-0402044 and DMS-$0701089$.

\section{\bf Results on functoriality which we need}

We have the
following theorem due to Jiang and Soudry, cf. Theorem 6.4 of
\cite{J-S1}, and Theorem 2.1 of \cite{J-S2}. As mentioned in the Introduction,
we are thankful to David Soudry for writing one of us (\cite{Sou-lett})
explaining the completion of their work; see also \cite{Sou} and \cite{GJR}.

\begin{thm}
Let $G$ be either a symplectic group or a quasi-split orthogonal group
over a non-Archimedean local field $k$.
Let the $L$-group of $G$
come equipped with its natural representation into $\GL_n({\Bbb C})$.
Then there exists a bijective correspondence between irreducible generic
discrete series
representations of $G(k)$ and irreducible generic representations
of $\GL_{n}(k)$, with Langlands parameter of the form
$$\sigma = \sum \sigma_i,$$
where $\sigma_i$ are pairwise inequivalent, irreducible representations of
$W_k^\prime \to {\rm GL}_n(k)$ which are all orthogonal except when
$G$ is an odd orthogonal group in which case they are all symplectic.

\end{thm}

We will also need following theorem due to Vign\'eras \cite{V}.

\begin{thm}

Let $K$ be a number field, and let $G$ be a quasi-split semi-simple group
over $K$. Let $v_i, i=1,\cdots, d$ be places of $K$, and $K_{v_i}$
the corresponding local fields. Suppose that $\pi_i$ are irreducible
generic supercuspidal  representations of $G(K_{v_i})$. Then there
exists a generic cuspidal automorphic representation $\Pi$ of
$G({\mathbb A}_K)$  with $\pi_i$  as the local component of $\Pi$ at
the place $v_i$, for $i=1,\cdots, d.$

\end{thm}

\medskip

The following theorem of Cogdell, Kim,
Piatetski-Shapiro, and Shahidi (cf. \cite{CoPSS}, \cite{CoKPSS})
establishes the weak Langlands functoriality from $G$ to GL$_n$.

\begin{thm}
Let $G$ be either a symplectic group or a quasi-split orthogonal group
over a number field $K$.
Let ${\mathcal A}^{0,g}(G({\mathbb A}_K))$ be the set of
irreducible generic cuspidal automorphic representations
of $G({\mathbb A}_K)$, and let ${\mathcal A}(\GL_{n}({\mathbb A}_K))$ be
the set of irreducible automorphic representations
of $\GL_{n}({\mathbb A}_K)$. Then there is a weak functorial lift
from
 ${\mathcal A}^{0,g}(G({\mathbb A}_K))$ to ${\mathcal A}(\GL_{n}({\mathbb A}_K))$.

\end{thm}

The weak Langlands functoriality established by Cogdell, Kim,
Piatetski-Shapiro, and Shahidi
has been proved to be Langlands functorial at all places of $K$ by
Jiang and Soudry, cf. Theorem E of \cite{J-S2} when $G={\rm SO}(2m+1)$,
and \cite{Sou}, \cite{Sou-lett}
for the remaining cases; a general statement is also in \cite{GJR}.

\begin{thm}
Let $G$ be either a symplectic group or a quasi-split orthogonal group
over a number field $K$.
The weak lift from  the set of irreducible generic
cuspidal automorphic representations
of $G({\mathbb A}_K)$ to the set
of irreducible automorphic representations
of $\GL_{n}({\mathbb A}_K)$ is Langlands functorial at every place of $K$,
i.e., it is the functorial lift at all
the places of $K$.

\end{thm}

An implication for us of these theorems is the following:

\begin{thm}

Let $K$ be a number field, and  $v_i, i=1,\cdots, d$  places of $K$.
Let $K_{v_i}$ be
the corresponding local fields.
Suppose that $\pi_i$ are irreducible
selfdual supercuspidal representations of $\GL_{n}(K_{v_i})$  whose
parameters are either  orthogonal, or symplectic for all $i$. Then there
exists a selfdual
cuspidal automorphic representation $\Pi$ on $\GL_{n}({\mathbb A}_K)$
with $\pi_i$ as the local component of $\Pi$ at the places $v_i$.

\end{thm}

\noindent{\bf Proof:} Since the parameters of the representations
 $\pi_i$ are all orthogonal or are all symplectic,
they can be transported, thanks to Theorem 1.1, to generic supercuspidal representations $\Pi_i$
of $G(K_{v_i})$ for a quasi-split orthogonal or symplectic group $G$ over a
number field $K$.
Using Theorem 1.2, these can be globalized into a
cuspidal automorphic representation $\Pi$ on $G(\A_K)$, i.e., with $\Pi_{v_i}\simeq \Pi_i$
for all $i$. The automorphic representation $\Pi$
can then be transferred to $\GL_{n}$ by Theorem 1.3, such that the lifted automorphic
representation $\pi$ on $\GL_{n}(\A_K)$ has the
correct local components $\pi_i$ at all the $v_i$ by Theorem 1.4,
proving Theorem 1.5.

\qed

\medskip

We note that such a globalization theorem will also follow from a stabilization of Arthur's
twisted trace formula for GL$(n)$. Some partial results are found in \cite{Ch-Cl}, which do not
suffice for us, however.

\medskip

\noindent{\bf Remark 1.6} \, It turns out that in our proof of
Theorem A, we will also need an analogue of Theorem 1.5 when one of
the $\pi_i$ is a (suitable) generalized Steinberg representation. A
soft result, sufficient for our purposes, is achieved in section 3
(see Theorem 3.1) using certain results of Clozel \cite{Cl2, Cl},
Harris \cite{HaLab}, Labesse \cite{Lab}, and Moeglin \cite{Moeglin}
involving the base change from unitary groups to GL$(n)$.

\medskip

\section{\bf Positivity of the global sign}

\medskip

Using the Jacquet-Langlands correspondence (see Theorem 2.5),
we can transport a selfdual, cuspidal
automorphic representation $\pi$ of $\GL_{n}({\mathbb A}_K)$ to a
selfdual automorphic representation $\pi^D$ of $D_K^\times({\mathbb A})$ of
discrete type, where $D_K$ is the central division algebra over $K$
ramified at only two places $v_1,v_2$ with invariants $\frac{1}{n}$,
and $-\frac{1}{n}$. The assertion of Theorem C is that $\pi^D$ is
orthogonal, i.e., that the global sign $c(\pi^D)$ is $1$.
We will prove this in a more general setup.

Let $G$ be a group, equipped with an involution $\theta$.
Call a representation $\eta$ of $G$ to be
$\theta$-selfdual if and only if $\eta^\vee\simeq\eta^\theta$. Here
$\eta^\theta$ is defined as $g\mapsto \eta(g^\theta)$, for all $g$ in
$G$. This implies that there exists a $G$-invariant bilinear
form $B': \eta \times \eta^{\theta}  \rightarrow {\mathbb C}$ which by
Schur's lemma is unique up to scaling. The bilinear form can also be
thought of as a bilinear form $B: \eta \times \eta \rightarrow {\mathbb
C}$ such that $B(gv,\theta(g)w) = B(v,w)$ for all $g \in G$, $v, w
\in \eta$. As $\theta$ is of order $2$, it is clear that $\tilde{B}$
defined by $\tilde{B} (v,w) = B(w,v)$ also has the same invariance
property with respect to $G$, and thus $\tilde{B} = cB$ for some
nonzero $c \in {\mathbb C}$. Since $\tilde{\tilde{B}} = B$, we find
that $c^2=1$, i.e., $c = \pm 1$. Thus we have associated an
invariant $c=c_\theta(\eta) \in \{ \pm 1 \}$ for any representation
$\eta$ of $G$ with $\eta^\vee \cong \eta^{\theta}$. The
analogue of the invariant $c_\theta(\pi)$ for Galois representations
appearing in the context of unitary groups was introduced by
Rogawski in chapter 15 of his book \cite{Rog}, see lemma 15.1.1, and
15.1.2. When $\theta = 1$, $c_\theta(\pi)$ supplies exactly the
information on whether $\pi$ is orthogonal or symplectic.

\begin{thm} Let $G$ be a reductive algebraic group over a number field $K$
together with an involution $\theta$.  Suppose $\Pi = \otimes_v^\prime
\Pi_v$ is an irreducible, $\theta$-selfdual automorphic
representation of $G(\A_K)$ of central character $\omega$, which
occurs discretely in $L^2(G(K)Z(\A_F)\backslash
G(\A_F),\omega)$ with multiplicity one. Then we have
$$
c_\theta(\Pi) \, = \, \prod\limits_{v\in {\rm ram}(G, \Pi)} \,
c_\theta(\Pi_v) \, = \, 1,
$$
where ram$(G, \Pi)$ denotes the finite set of places outside which
$G_v$ is quasi-split  and $\Pi_v$ is unramified.
\end{thm}

\medskip

\noindent{\bf Proof} \, Define a $G({\mathbb A})$-invariant bilinear
form $B$ on $\Pi$ by
$$ (f,g) \longrightarrow
\int_{G(K)Z(\A_F)\backslash G(\A_F)}fg^\theta d \mu ,
\,\,\,\,\,\,\,\,\,\, \forall \, f,g \in \Pi,
$$
where $d\mu$ is an invariant measure on $G(K)Z(\A_F)\backslash
G(\A_F)$. We check that this is a non-degenerate bilinear form on
$\Pi$. Note that the space of functions spanned by $\bar{f}$ (the
complex conjugate of $f$), as $f$ varies over $\Pi$, gives rise to
the representation $\Pi^\vee$, which is isomorphic to $\Pi^\theta$.
Hence by the multiplicity 1 hypothesis, $\bar{f}^\theta \in \Pi$.
Taking $g= \bar{f}^\theta$, we see that
$$
B(f,g) \, = \, \int f \bar{f} d \mu \, \ne \, 0.
$$
Consequently, the bilinear form $B$ on $\Pi$ is non-degenerate. We
have to show that $c_\theta(\Pi)=1$. (When $\theta=1$, this is
obvious as the bilinear form is evidently symmetric.) It suffices to
show that there exist $f, g$ such that $B(f,g)$ and $B(g,f)$ are
both positive. But if we take $g=\bar{f}^\theta$ as above, then
$B(f,g)$ (resp. $B(g,f)$) is the $L^2$-norm of $f$ (resp.
$f^\theta$). So the $\theta$-sign of $\Pi$ is positive.

It is left to show that for almost all places $v$,
$c_\theta(\Pi_v)=1$. This is a consequence of the following result,
observed in \cite{PR} for $G=\GL_n$.

\begin{prop} Every irreducible admissible
representation $\pi$ of a quasi-split reductive group $H$
over a local field, satisfying $\pi^\vee\simeq\pi^\theta$ for an
involution $\theta$ of $H(k)$, has $\theta$-sign $1$ whenever $\pi$
is unramified (for a maximal compact subgroup left invariant by $\theta$).
When $H=$\GL$(n)$ and $\theta=1$, the same assertion
holds for any $\pi$ (not necessarily generic, by the theory of
degenerate Whittaker models).

\end{prop}

Theorem 2.1 now follows, as does Theorem C.

For our applications, we will need some information at infinity, where
we will assume, and this suffices for the purposes of this paper,
 that the group is $\GL_n({\Bbb R})$. We note the following
lemma.

\begin{lemma}
Let $\pi$ be an irreducible $(\mathfrak{g},K)$-module
where $\mathfrak{g}$ is the Lie algebra of $\GL_n({\R})$ (resp.
$\GL_n({\C})$), and $K =\O_n({\Bbb R})$ (resp. $U_n$)
is the maximal compact subgroup of $\GL_n({\R})$ (resp.
$\GL_n({\C})$). Then if $\pi$
is selfdual, it carries a symmetric bilinear form, i.e., $c(\pi)=1$.
\end{lemma}

\noindent{\bf Proof:} It is a general result due to D.~Vogan that any
irreducible $(\mathfrak{g},K)$-module has a minimal $K$-type
which occurs with multiplicity $1$.
(All that matters for us is that there is a $K$-type which
 appears with multiplicity 1.) Since
it is well-known that every irreducible, selfdual representation of $\O_n({\R})$ (resp. $U_n$)
carries an invariant symmetric bilinear form,
the lemma follows. (Note the difference between $U_n$ and $SU_n$; there are already irreducible
symplectic representations of $SU_2$, but they do not extend to selfdual representations of $U_2$.)

\bigskip

\section{\bf An interlude on globalization}

\medskip

It appears difficult to construct, for general $n$, selfdual cusp forms on GL$(n)$ which have a given generalized Steinberg component at one place and a supercuspidal component at another. The following represents a stopgap measure, which happens to be sufficient for our purposes (see section 4).

\medskip

\begin{thm}\label{global} Fix $n=4dm$, with $d, m \geq 1$, and let $k$ be a non-archimedean local field. Let $\pi$ be a discrete series representation of GL$_n(k)$ with Langlands parameter $\sigma = \sigma_0\otimes sp_{2m}$, where $\sigma_0$ is irreducible and induced by a character $\chi$ of the unramified degree $2d$ extension $k_{2d}$ of $k$ such that $\chi^\alpha=\chi^{-1}$, with $\alpha$ denoting the unique element of order $2$ of Gal$(k_{2d}/k)$. Then there exists a totally real number field $F$ with $F_v=k$ for a finite place $v$, and a selfdual, cuspidal automorphic representation $\Pi$ of GL$_n(\A_F)$, such that (i) $\Pi_v\simeq \pi$ and (ii) $\Pi_u$ is, for another finite place $u$ of $F$, a supercuspidal representation $\pi'$, say, with parameter $\sigma'$ satisfying $c(\sigma')=c(\sigma)$.
\end{thm}

{\it Proof of Theorem 3.1}. \, Since $sp_{2m}$ is the symmetric $(2m-1)$-th power of the standard representation of SL$(2,\C)$, it is symplectic, and thus
$$
c(\sigma) \, = \, c(\sigma_0)c(sp_{2m}) \, = \, -c(\sigma_0).
$$
We begin with a simple lemma.
\begin{lemma} Fix $m\geq 1$. Let $k'$ be a non-archimedean local field of residual characteristic $p$ odd, which contains all the $2m$-th roots of unity. Then there exists an irreducible $2m$-dimensional, symplectic representation $\tau$ of $W_{k'}$ such that, for any finite unramified extension $k''/k'$, the restriction of $\tau$, to $W_{k''}$ remains irreducible.
\end{lemma}

{\it Proof of Lemma}. \, Let $E/k'$ be a ramified cyclic $2m$-extension,
which exists by Kummer theory, and write $L$ for the unique cyclic
$m$-extension of $k'$ contained in $E$. It suffices  to construct a character
$\lambda$ of $E^\times$ such that its restriction to $L^\times$ is the
quadratic character $\omega_{E/L}$ associated to the quadratic extension $E$
by the class field theory, and such that its restriction to the units
${\mathcal O}_E^\times$ has $2m$ distinct Galois conjugates. For this observe
that the exponential mapping induces an isomorphism, which is Galois
equivariant, between the additive group ${\mathcal P}_E^r$ and the
multiplicative group $1+{\mathcal P}_E^r$ for sufficiently large
integer $r$, and therefore to construct a character on $1+{\mathcal P}_E^r$,
it suffices to construct one on the additive group ${\mathcal P}_E^r$.
Now observe that by the normal basis theorem, $E$ is a free $k'[G]$-module
of rank 1 where $G$ is the Galois group of $E$ over $k'$. From this
it is easy to construct a homomorphism of ${\mathcal O}_{k'}$-modules
from  ${\mathcal P}_E^r$ to ${\mathcal O}_{k'}$ which is trivial
on ${\mathcal P}_L^r$, and has $2m$ distinct Galois conjugates. Composing
with an appropriate character of ${\mathcal O}_{k'}$ to ${\Bbb C}^\times$,
we obtain a character on $1+{\mathcal P}_E^r$ which has $2m$
Galois conjugates, and which is trivial on $1+{\mathcal P}_L^r$. One
can extend this character to $E^\times$ such that its restriction to
$L^\times$ is $\omega_{E/L}$, giving us the desired character $\lambda$. We
set $\tau:={\rm Ind}_E^{k'}(\lambda)$.

\qed

Now we choose a totally real number field $F$ with finite places $u, v$ such that $F_v=k$ and $F_u=k'$. Then we pick a cyclic extension $M/F$ of degree $2d$ in which both $v$ and $u$ are inert (of degree $2d$ over $F$), such that the unique quadratic extension $K$ of $F$ contained in $M$ is a CM field. Denote the unique prime divisors of $v$ in $M$, resp. $K$,  by $\tilde w$, resp. $w$, and the corresponding ones over $u$ by $\tilde w'$, resp. $w'$. By construction, $M_{\tilde w}=k_{2d}$, $K_{w}=k_2$, $M_{\tilde w'}=k'_{2d}$, and $K_{w'}=k'_2$.

Let $\eta$ denote the supercuspidal representation of GL$_{2m}(k')$ associated to $\tau$ by the local Langlands conjecture. Fix a conjugate selfdual character $\chi'$ of ${k'_{2d}}^\times$ whose restriction to ${k'_d}^\times$
is trivial if and only if the restriction of the given character $\chi$ of $k_{2d}^\times$ is trivial on $k_d^\times$. This implies that $c(\sigma_0)=c(\sigma_0')$, where
$$
\sigma_0'={\rm Ind}_{k'_{2d}}^{k'}(\chi').
$$
We put
$$
\sigma' \, \simeq \, \sigma_0' \otimes \tau,
$$
and note that it is irreducible by the Lemma. Since by construction, $c(sp_{2m})=c(\tau)=-1$, we also have
$$
c(\sigma) \, = \, c(\sigma').
$$
Now let $\pi'$ be the supercuspidal representation of GL$_n(k')$, $n=4dm$, which is associated to $\sigma'$ by the local Langlands correspondence. Note that
$$
\pi' \, \simeq \, I_{k'_{2d}}^{k'}(\chi')\boxtimes\eta,
$$
where $\boxtimes$ denotes the functorial product.
It is this $\pi'$ for which Theorem 3.1 holds.

\medskip

We now need the following two auxiliary results, which are likely known to experts, but we indicate their proofs for completeness.

\begin{thm}\label{grossen} Let $M/N$ be a quadratic extension of number fields with non-trivial automorphism $\alpha$.
Denote by $S$ a non-empty
finite set of finite places $x$ of $M$ (lying above places $y$ of $N$) such that the local degrees $[M_x:N_y]$ are all $2$.
For every $x$ in $S$, let
$\chi_x : M_x^\times \rightarrow {\C}^\times$
be a finite order character such that $\chi_x^\alpha = \chi_x^{-1}$.
Then $\{\chi_x \, \vert \, x\in S\}$ can be globalized to a
$\alpha$-self-dual character $\chi$ of
${\A}^\times_M/M^\times$ if and only if
the restrictions $\chi_x\vert_{_{{N_y^\times}}}, \, x\in S$, are all trivial
or all non-trivial.

\end{thm}

\begin{thm}\label{conjugate global} Let $K$ be a CM quadratic extension of a totally real number field $F$, with non-trivial automorphism $\theta$ over $F$. Fix places $v, u$ of $F$ which are inert in $K$ with respective divisors $w, w'$; put $k=F_v$ and $k'=F_u$. Let $St$ denote the Steinberg representation of GL$_{2m}(K_w)$ and $\eta$ a selfdual, supercuspidal representation of GL$_{2m}(k')$ whose parameter $\tau$ is symplectic. Then there exists a $\theta$-selfdual, cuspidal automorphic representation $\Pi_1$ of GL$_{2m}(\A_K)$ such that $\Pi_{1,w}$ is $St$ and  $\Pi_{1,w'}$ is $\eta_2$, the base change of $\eta$ to $K_{w'}=k'_2$.
\end{thm}

\medskip

\noindent{\bf Claim 3.5}: \, {\it Theorem 3.3 $+$ Theorem 3.4 \, $\implies$ \, Theorem 3.1}

\medskip

{\it Proof of Claim 3.5}. \, Preserve the hypotheses of Theorem 3.1, and the constructions in the paragraphs following its statement. Globalize $\chi$ and $\chi'$ simultaneously, using Theorem 3.3, to an $\alpha$-selfdual character $\Psi$ of $M$ such that $\Psi_{\tilde w}=\chi$ and $\Psi_{\tilde w'}=\chi'$. Next apply Theorem 3.3 to deduce the existence of a $\theta$-selfdual, cuspidal automorphic representation $\Pi_1$ of GL$_m(\A_K)$ such that $\Pi_{1,w}\simeq St$ and $\Pi_{1,w'}\simeq \eta_2$. Put
$$
\Pi : = \, I_M^F(\Pi_{1, M}\otimes\Psi) \, \simeq \, I_K^F(\Pi_1\boxtimes I_M^K(\Psi)),
$$
where $I_M^F$ (resp. $I_K^F$) denotes the (global) automorphic induction from GL$(2m)/M$ (resp. GL$(2dm)/K$) to GL$(n)/F$ (\cite{AC}); see also \cite{HH}), and the functorial product $\boxtimes$ makes sense here because $I_M^K(\Psi)$ is cyclic monomial. By construction, $\Pi_v \simeq \pi$ and $\Pi_u \simeq \pi'$. Moreover, since $\pi'$ is supercuspidal, $\Pi$ must be cuspidal. Finally, since $\Pi_1$ and $I_M^K(\Psi)$ are both $\theta$-selfdual, $\Pi$ is forced to be selfdual.

\qed

\medskip

{\it Proof of Theorem 3.3}. \, The condition $\chi^\alpha = \chi^{-1}$ implies that the
character $\chi$ must be, upon restriction to ${\A}^\times_N/N^\times$, trivial
on the index $2$ subgroup consisting of norms from ${\A}^\times_M/M^\times$,
implying that $\chi|_{\A^\times_N}$ is either 1, or is equal to the quadratic character
$\omega=\omega_{M/N}$ of ${\A_N}^\times$ associated to the extension $M/N$. Thus the
existence of a global, $\theta$-selfdual character $\chi$ implies that
$\chi|_{N_y^\times}$ is either the trivial character
$\forall x \in S$, or is the non-trivial quadratic character of $\omega_y$
(associated to the extension $M_x/N_y$) $\forall x \in S$.

Conversely, assume that $\chi_x|_{N_y^\times} \equiv 1$ for all $x\vert y$ in $S$. Then
by Hilbert Theorem $90$, there are characters $\eta_x$ of $M_x^\times$, $x \in S$,
such that $\chi_x
= \eta_x^\theta\eta^{-1}$. By the Grunwald-Wang
theorem, these characters $\eta_x$, $x \in S$, can be globalized to a
finite order character $\eta$ of $\A_M^\times/M^\times$. If we put $\chi:=\eta^\theta\eta^{-1}$.
then it has the given local components $\chi_x$ at all the $x$ in $S$, and moreover, it is $\theta$-selfdual.

Assume next that $\chi_x|_{N_y^\times} =\omega_y$ for all $x\vert y$ in $S$.
Let $\tilde{\omega}$ be a character of ${\A}^\times_M/M^\times$
which restricts to $\omega$ on ${\A}^\times_N/N^\times$. For every $x$ in $S$,
put $\psi_x = \chi_x\tilde\omega_x^{-1}$, which restricts to the trivial character on $N_y^\times$
for all $x\in S$. So, by the previous case, there is a $\theta$-selfdual character $\Psi$ of
${\A}^\times_M/M^\times$ with components $\psi_x$ at $x \in S$. We are now done by setting $\chi=\Psi\tilde\omega^{-1}$.

\qed

\medskip

{\it Proof of Theorem 3.4}. \, Preserving the hypotheses of the Theorem, denote by $\tilde \omega$ an idele class character of the CM field $K$ whose restriction to the idele class group of $F$ is the quadratic character $\omega=\omega_{K/F}$ attached to $K/F$; then $\tilde\omega$ is necessarily conjugate selfdual. Fix a finite set $S$ of finite places of $K$ which are of degree $1$ over $F$ with $\vert S\vert \geq 2$. Since $\eta$ is a supercuspidal representation of GL$_{2m}(F_u)$ of symplectic type, we may, by using Theorem 1.5,  globalize it to be the $u$-component of a selfdual, cuspidal automorphic representation $\Pi_0$ of GL$_{2m}(\A_F)$, such that the set of places where $\Pi_0$ is supercuspidal includes, besides $u$, places of $F$ below $S$. Then the base change $\Pi_{0,K}$ to GL$(2m)/K$ is still cuspidal as it has supercuspidal components, and moreover, it is $\theta$-selfdual.

Let $G$ denote an anisotropic unitary group  over $F$ associated to $K/F$ such that $G_v$ and $G_w$ are quasi-split; it may be defined either by a hermitian form in $2m$ variables over $K$ or by a division algebra $\D/K$
of dimension $m^2$, ramified only at a subset of $S$ and split at all the other places (including $w, w'$), and equipped with an involution $\tilde\theta$ which restricts to $\theta$ on $K$. Then by \cite{Cl, Lab} or \cite{HaLab}, Theorem 2.4.1, we know that either $\Pi_{0,K}$ or $\Pi_{0,K}\otimes \tilde\omega$ descends to a cuspidal automorphic representation $\pi_0^G$ of $G(\A_F)$, which is functorial almost everywhere and is compatible with local Langlands transfer where it is available. Since $\Pi_{0,w'}\simeq\eta_2$ is a discrete series, we know by Moeglin \cite{Moeglin}, cf. sections 6-8,  that the local descent $\beta:=\pi^G_{0,u}$ belongs to a stable discrete series $L$-packet of $G(F_u)$ attached to the parameter $\psi: W'_{k'}\to {}^LG$ defined by $\Pi_{0,w'}$. (In \cite{Moeglin}, she considers a larger class of ''$\theta$-discrete'' representations, and associates canonically an elliptic, stable $L$-packet on a unique twisted endoscopic group of the form $U(n_1)\times U(n_2)$, but in the discrete series case either $n_2$ or $n_1$ is zero, and which case occurs depends on whether the representation itself, or its $\tilde\omega_{w'}$-twist, descends.) Since $\Pi_{0,w'}$ is conjugate selfdual, its parameter defines an admissible homomorphism $\psi_2: W'_{k'_2}\to {}^LG$ whose image is contained in no proper Levi, and either $\psi_2$ or $\psi_2\otimes\tilde\omega_{w'}$, but not both, extends to give the parameter $\psi$ on $W'_{k'}$. Let us write $\mu=1$ or $\tilde\omega$, so that $\Pi_{0,K}\otimes\mu$ descends to $\pi_0^G$. By Lemma 15.1.2 of \cite{Rog}, $\mu_{w'}$ is trivial on ${k'}^\times$ if and only if the $\theta$-sign of its parameter, namely $c_\theta(\tau_2)$, with $\tau_2$ denoting $\tau_{\vert_{k'}}$, is $1$. There is also a weak base change (see \cite{Cl} or \cite{HaLab}, Theorem 2.2.2), transferring $\pi_0^G$ back to a cusp form $\Pi'$ on GL$_{2m}(\A_K)$, which is conjugate selfdual and is functorial at almost all places. (We use Jacquet-Langlands in addition to base change if $G$ is defined by $\D$.) Then by the strong multiplicity one theorem for GL$(2m)/K$, $\Pi'$ and $\Pi_{0,K}\otimes\mu$ must be isomorphic, and in particular, $\Pi'_{w'}$ is just the supercuspidal representation $\eta_2\otimes\mu_{w'}$.

Next we take up the problem of globalizing $\eta_2$ and $St$ simultaneously to a conjugate selfdual cusp form $\Pi_1$ on GL$(2m)/K$. We construct $\Pi_1\otimes\mu$ as the base change of a suitable cusp form $\pi^G$ on $G/F$. Note that by hypothesis, $c(\tau)=c(sp_{2m})=-1$, which implies that $c_\theta(\tau_2)=c_\theta(sp_{2m})$, and so we have no obstruction, and $St\otimes\mu_w$ will descend to $G(F_v)$.

\medskip

\noindent{\bf Proposition 3.6} \, \it There exists a cuspidal automorphic representation $\pi^G$ of $G(\A_F)$ such that its components at $v, u$ and the archimedean places are sufficiently regular discrete series representations, and moreover, up to modifying $\mu$ by a character which restricts to the trivial character of $\A_F^\times/F^\times$, the local component $\pi^G_v$, resp. $\pi^G_u$, base changes to $St\otimes\mu_w$, resp. $\eta_2\otimes\mu_{w'}$, on ${\rm GL}(2m)/K_w$, resp. on ${\rm GL}(2m)/K_{w'}$.
\rm

\medskip

{\it Proof of Proposition 3.6}. \, Let $St^G$ denote the Steinberg representation of $G(F_v)$. By appealing to a result of Clozel on limit multiplicities (cf. Theorem 2 of \cite{Cl2}), we get the existence of a cusp form $\pi^G$ with sufficiently regular discrete series (resp. supercuspidal) components at infinity (resp. $S$), such that $\pi^G_v\simeq St^G$ and $\pi^G_u\simeq \beta\in \Pi_u(\psi)$. Denote by $\Pi_1$ the base change of $\pi^G$ to GL$_{2m}/K$. As we have seen above, the local base change at $u$ sends $\beta$ to $\eta_2\otimes\mu_{w'}$. Moreover, one sees (cf. section 3.9 of \cite{Lab}) that the local base change $\Pi_{1,w}$ of $St^G$ to GL$_{2m}(K_{w}$ is an unramified twist of $St$. Indeed, the global base change is achieved by comparing the trace of functions $f$ on $G(\A_F)$ with the twisted trace of functions $\varphi$ on $\D^\times(\A_K)$, and locally at $v$, Proposition 3.9.2 of {\it loc. cit.} matches the Kottwitz functions $f_v$ on $G(F_v)$ and $\varphi_w$ on GL$_{2m}(K_w)$, and these functions are stabilizing and have zero traces on all the irreducible unitary representations except for those occurring as subquotients, up to unramified twists, of principal series representations admitting the trivial representation as the Langlands quotient. Note that unramified twists of non-Steinberg representations of GL$_{2m}(K_{w'})$ occurring as such subquotients are non-generic and cannot occur as local components of cusp forms on GL$(2m)/K$. It follows that $\Pi_{1,w}\simeq St\otimes\nu$ for an unramified character $\nu$, which must differ from $\mu_w$ by a character which is trivial on $F_v^\times$.

\qed

\medskip

To conclude, we get a $\theta$-selfdual, cuspidal automorphic representation $\Pi_1$ of GL$_m(\A_K)$ associated to $\pi^D$, which, when twisted by $\mu^{-1}$, will satisfy all the assertions of Theorem 3.4.

\qed

\bigskip

\section{\bf Sign in the Level 1 case}

\medskip

In this section we prove Theorem A about
irreducible selfdual representations of $D^\times/D^\times(1)$. The global proofs in this
paper depend most crucially on this local input, as it represents the simplest of the situations
for Theorem A.

Various aspects of the representation theory of $D^\times/D^\times(1)$ are analyzed in the
work of Silberger and Zink in \cite{S-Z}. We begin with  some notation, and recalling
the parametrization of the irreducible representations of $D^\times/D^\times(1)$ from \cite{S-Z},
where they unfortunately call them level zero representations.

As in the rest of the paper, let $D$ be a division algebra with center a
non-Archimedean local
field $k$, and of index $n$. Let $n = e f$, and let $k_f$
be the unramified extension of $k$ of degree $f$, contained in $D$. Let $D_f$ be the
centralizer of $k_f$ in $D$ which is a division algebra with center $k_f$ and
of index $e$. A character $\chi$ of $k_f^\times$ will be called regular if all its
Galois conjugates are distinct. For a character $\chi$ of $k_f^\times$, let $\tilde{\chi}$
be the character of $D_f^\times$ obtained by composing with the reduced norm mapping
$ {\rm Nrd} : D_f^\times \rightarrow k_f^\times$. If the character $\chi$ is tame, i.e., trivial
on $k_f^\times(1)$, then the character $\tilde{\chi}$ of $D_f^\times$ can be extended to a
character of $D^\times(1)D_f^\times$ by declaring it to be trivial on $D^\times(1)$, which, by
abuse of notation, we again denote by $\tilde{\chi}$.

With this notation, it follows from Clifford theory as in \cite{S-Z} that irreducible
representations of $D^\times/D^\times(1)$ have dimensions $f$ which are divisors of $n$, and
that there is a bijection between irreducible representations of $D^\times/D^\times(1)$
of dimension $f$ and regular characters of $k_f^\times$ which are trivial on $k_f^\times(1)$
(modulo the action of the Galois group of $k_f$ over $k$ on the set of such characters) obtained by inducing
the character $\tilde{\chi}$ of $D^\times(1)D_f^\times$ to $D^\times$.

To analyze whether these representations, when selfdual, are orthogonal or symplectic,
we state the following general proposition, whose straightforward proof will be omitted.

\begin{prop}
Let $H$ be a normal subgroup of a group $G$ of index $f>1$, such that $G/H$ is a cyclic group
of order $f$. Let $\varpi$ be an element of $G$ whose image in $G/H$ is a generator
of the cyclic group $G/H$.  Let $\pi$ be an irreducible   representation of $G$ of dimension
$f$ whose restriction to $H$ contains a character $\chi: H \rightarrow {\Bbb C}^\times$, so that
$\pi = {\rm Ind}_H^G \chi$. Then,

\begin{enumerate}
\item If $\pi$ is selfdual, $f$ is even, say $f=2d$, which we assume is the case
in the rest of the proposition.

\item The representation $\pi$ is selfdual if and only if
 $\chi^{-1}= \chi^{<d>}$ where $\chi^{<d>}(h)= \chi(\varpi^d h
\varpi^{-d})$ for $h \in H$.

\item $\varpi^{f}$ acts on $\pi$ by a scalar $C(\pi)$ on the representation space $\pi$. If
$\pi$ is selfdual, $C(\pi) = \pm 1$, and $C(\pi) = 1$ if and only if $\pi$ is an orthogonal
representation.

\item If $\pi$ is selfdual, it is orthogonal if and only if $\det \pi(\varpi) = -1$.
\end{enumerate}
\end{prop}

Continuing now with the representations of $D^\times$ of dimension $f$ with $ef $, let $\varpi$
be an element of $D^\times$ which normalizes $k_n$, an unramified extension of $k$ inside $D$
of degree $n$ over $k$, such that $\varpi^n = \varpi_k$, a uniformizer in $k$.The element
$\varpi$ of $D^\times$ projects to a generator
of the cyclic group $D^\times/D^\times(1)D_f^\times \cong {\Bbb Z}/f{\Bbb Z}$, and as $\varpi^f$
centralizes $k_f$, it lies in $D_f$. Since $(\varpi^f)^e=\varpi_k$, it follows that the reduced
norm of $\varpi^f$ is $\varpi_k$. From the previous proposition, we conclude the following corollary.

\begin{cor}
Let $\chi$ be a regular character of $k_f^\times$, and $\pi_\chi$ the associated representation
of $D^\times$ of dimension $f$.
Then $\pi_\chi$ is selfdual if and only if $f$ is even, say $f=2d$, and the character
$\chi$ restricted to $k_d^\times$ is trivial on the index 2 subgroup consisting of norms
from $k_f^\times$.  Assuming $\pi_\chi$ to be selfdual, it is orthogonal if and only if
$\chi$ restricted to $k^\times$ is trivial.
\end{cor}

Recalling that the Weil group $W_{k_f/k}$ sits in the exact sequence,
$$1 \rightarrow k_f^\times \rightarrow W_{k_f/k} \rightarrow {\rm Gal}(k_f/k) \rightarrow 1,$$
and there is again an element $\varpi$ in $W_{k_f/k}$ which goes to the generator
of the Galois group of $k_f$ over $k$, and whose $f$-th power is a uniformizer in $k_f$,
we have the following corollary for representations of the Weil group.

\begin{cor} For a regular character $\mu$ of $k_f^\times$, let $\sigma_\mu$ be the induced
representation of $W_{k_f/k}$ of dimension $f$.  Then $\sigma_\mu$ is
selfdual if and only if $f$ is even, say $f=2d$, and the character
$\mu$ restricted to $k_d^\times$ is trivial on the index 2 subgroup consisting of norms
from $k_f^\times$.  Assuming $\sigma_\mu$ to be selfdual, it is orthogonal if and only if
$\mu$ restricted to $k^\times$ is trivial, or if and only if $\det \sigma_\mu$ is nontrivial.
\end{cor}

We will give a proof of Theorem A below without knowing the precise Langlands parameter
attached to $\pi_\chi$. In any case, we believe the following question has an affirmative answer:

\medskip

\noindent{\bf Question 4.4} \,
Does the Langlands parameter of the representation $\pi_\chi$ of dimension $f$ of $D^\times$
equal $\sigma_\mu \otimes sp_e$, with $\sigma_\mu$ the $f$-dimensional representation of $W_k$
induced by the character
$$
\mu:=\chi \omega_2^{e(f-1)}: k_f^\times\rightarrow \C^\times,
$$
where $\omega_2$ is the quadratic unramified character of $k_f^\times$?

\medskip

\noindent{\bf Remark 4.5 :} We have not been able to find a precise reference dealing with this question.
However, we should point out that in \cite{S-Z} Silberger and Zink suggest without proof,
see the Remark on page 182 of {\it loc. cit.},
that the Langlands parameter of $\pi_\chi$ is $\sigma_\mu\otimes sp_e$,
with $\mu = \chi \omega_2^{f-1}$. There is no discrepancy with our formula above unless $e$ and $f$ are both even. Results of 4.6 below show, however, that their formula cannot be correct for all $e$. Finally, Henniart has informed us that he is working on a paper with Bushnell to answer Question 4.4.

\medskip

\subsection*{\bf 4.6 \, Proof of Theorem A for level $1$ representations}

\medskip

Let $\pi$ be an irreducible representation of $D^\times$ of level $1$ as above. The parameter
of $\pi$ is necessarily of the form
$$
\sigma =\sigma_\nu \otimes sp_e , \quad {\rm with} \quad \sigma_\nu= {\rm Ind}_{W_{k_{2d}}}^{W_k}(\nu)\otimes sp_e,
$$
for a tame character $\nu$ of $W_{k_{2d}}$ satisfying $\nu^\theta=\nu^{-1}$, where $\theta$ denotes the non-trivial automorphism of $k_{2d}/k_d$. The irreducibility of $\sigma$ also forces $\nu$ to be distinct from $\nu^\xi$ for
any $\xi$ in Gal$(k_{2d}/k)$.

First suppose $e$ is $1$. Since the central character of $\pi_\chi$ is associated, by the local correspondence, to the determinant of $\sigma_\nu$, the former is trivial if and only if the latter is also trivial. The Corollaries 3.2 and 3.3 then imply that $\pi_\chi$ is orthogonal if and only if $\sigma_\nu$ is symplectic. This argument works for {\it $e$ odd} in the same way.

\medskip

We may now assume that $e=2m$ and $f=2d$ are both {\it even}, so that $n=4dm$. Thanks to Theorem 3.1, we can find a totally real number field $F$ with $F_v=k$ at a place $v$, and a cuspidal, selfdual automorphic representation $\Pi$ of GL$_n(\A_F)$ with  $\Pi_v$ in the discrete series corresponding to $\pi$, and at another finite place $u$, $\pi':=\Pi_u$ is supercuspidal. Let $\B/F$ be a division algebra of index $n$ over $F$ which ramifies only at $u, v$ such that $\B_v=D$ and $\B_u=D'$, where $D'$ is also a division algebra. Applying Theorem 2.1 to the automorphic representation $\Pi^\B$, say, of $\B(\A_{F_0}^\times)$ associated to $\Pi$, we conclude, since $\beta_{v_0}\simeq\pi$, that $c(\pi)=c(\pi')$, where $\pi'$ corresponds to the supercuspidal representation $\Pi_u$. By using Theorem 1.5, we can also find a cuspidal, selfdual automorphic representation $\Pi_0$ of GL$_{2m}(\A_F)$ such that $\Pi_{0,u}=\Pi_u$ and $\pi'':=\Pi_{0,v}$ supercuspidal, with the parameters of these two supercuspidals being both symplectic or both orthogonal. In fact, we have freedom in choosing $\pi''$, and we may arrange it to correspond to a level $1$ representation of $D^\times$. Applying Theorem 2.1 again, we get $c(\pi'')=c(\pi')$, and since the parameter $\sigma''$ of $\pi''$ is trivial on SL$(2,\C)$, we know by what we did above that $c(\pi'')=-c(\sigma'')$. It follows that $c(\pi)$ equals $-c(\sigma)$ as desired.

\qed

\noindent{\bf Remark 4.7}: \, To be able to make use of representations of
level $1$,
we must know that there is an irreducible selfdual representation
of $D^\times/D^\times(1)$ of dimension $f$ for any {\it even}
divisor $f=2d$ of $n$. This reduces to a question over finite fields.
A character of ${\Bbb F}_{q^{2d}}$ is regular if and only if it does not
factor through the norm mapping of an intermediate field, say ${\Bbb F}_{q^s}$, in which case it has an order divisible by $(q^s-1)$, whereas a
character of ${\Bbb F}_{q^{2d}}$ gives rise to a selfdual representation
if and only if  it arises from
the circle group of norm 1 elements of ${\Bbb F}_{q^{2d}}$, denoted
${\Bbb S}^1({\Bbb F}_{q^d})$, through the map $x\rightarrow x/\bar{x}$ for
$x \in {\Bbb F}_{q^{2d}}$.
Since ${\Bbb S}^1({\Bbb F}_{q^d})$ is a cyclic group of order $(q^d+1)$, one
can consider characters on it of order $(q^d+1)$, which will then not arise
from an intermediate field through the norm mapping.
\bigskip

\section{\bf Proof of Theorem A for $D^\times$}

\medskip

Let $\pi'$ be an irreducible, selfdual representation of $D^\times$ of
parameter $\sigma$, and associated discrete series representation $\pi$
of GL$_n(k)$. First assume that $\pi$ is supercuspidal.

When $n$ is {\it odd}, embed $\pi$, by using Theorem 1.5, as a local component $\Pi_v$ of a
selfdual, cuspidal automorphic representation $\Pi$ of GL$_n(\A_K)$, for a number field $K$
with $K_v=k$. Let $E/K$ be a cyclic extension of degree $n$ such that $v$ splits completely
in $E$, and denote by $\Pi_E$ the base change of $\Pi$ to GL$_n/E$ as defined by Arthur and
Clozel (\cite{AC}). Let $\B$ be the division algebra of index $n$ over $E$, ramified only at the
places of $E$ above $v$, such that $\B_v:= \B\otimes_K K_v \simeq D^n$. (By class field theory,
such a division algebra exists.) Let $\Pi^{\B}$ be the automorphic representation of
$\B^\times(\A_E)$ associated to $\Pi_E$ by the Jacquet-Langlands correspondence (\cite{bad}),
which appears with multiplicity one in the space of automorphic forms on $\B^\times/E$. Then
the component of $\Pi^{\B}$ at $E_v:=E\otimes_K K_v \simeq K_v^n$
is necessarily
isomorphic to $\pi'^{\otimes n}$. Since at any place
$w$ of $E$ not lying over $v$, $\B_w^\times$ is by
construction GL$_n(E_w)$, we see by applying Theorem 2.1 that
$$
c(\pi')^n \, = \, c(\pi'^{\otimes n}) \, = \, c(\Pi^{\B}) \, = \, 1.
$$
It follows, by the oddness of $n$, that $c(\pi')=1$, proving
Theorem A in this case.

\medskip

Now suppose $n$ is even, still with $\pi$ being supercuspidal. Now globalize $\pi$, again using
Theorem 1.5, to a selfdual, cuspidal automorphic representation $\Pi$ of GL$_n(\A_K)$, with $K_v=k$
and $\Pi_v=\pi$. Moreover we can arrange this $\Pi$ in such a way that, at a second finite place $u$,
the local component $\Pi_u$ is a selfdual, supercuspidal representation of level $1$, i.e., corresponding to a
representation of {\it level $1$} of a division algebra
over $K_u$ of index $n$.
Now choose a global division algebra $\D$ of index $n$ over $K$ such that $\D$
is ramified only at $u, v$, with $\D_v\simeq D$, and denote by $\Pi^\D$ the automorphic representation of
$\D^\times(\A_K)$ associated to $\Pi$. Applying Theorem 2.1, we obtain
$$
c(\pi')c(\Pi^\D_u) \, = \, c(\Pi^\D) \, = \, 1.
$$
In section 4, we proved Theorem A in the level $1$ situation. Thus we get
$$
c(\pi') \, = \, c(\Pi^\D_u) \, = \, -c(\sigma_u) \, = \, -c(\sigma),
$$
as asserted.

We have now
achieved a proof of Theorem A for those representations $\pi'$ of $D^\times$
for which the corresponding representation $\pi$ of $\GL_n(k)$
is supercuspidal, thus the Langlands parameter is trivial on
the $\SL(2,{\mathbb C})$ part of $W_k^\prime$.
This restriction has been imposed on us as we do
not have the globalization theorem (Theorem 1.5) available for general
discrete series representations.

\medskip

Let $\pi$ be a non-supercuspidal discrete series representation of
$\GL_n(k)$. Then its parameter $\sigma$ must be of the form $\tau\otimes sp_b$,
$n = ab$, with $sp_b$ being the unique $b$-dimensional irreducible
representation of $\SL_2({\mathbb C})$ and $\tau$ an irreducible
selfdual $a$-dimensional representation of $W_k$.

Let $\Sigma$ be a selfdual cuspidal automorphic representation of
$\GL_a({\mathbb A}_K)$ whose local component at the place $v$ of $K$
with completion $k$ has Langlands parameter $\tau$. We may assume, thanks to
Theorem 1.5, that
at some other finite place, say $u$, $\Sigma_u$ is supercuspidal of level $1$,
whose parameter $\tau_{u}$ is of the same parity as $\tau$. We may take $\tau_u$ of the form
${\rm Ind}_{L_a^\times}^{W_k} (\chi_1) $, where
$\chi_1$ is a certain character of $L_a^\times$ trivial on $L_a^\times(1)$ where $L_a$ is the unique unramified
extension of $K_u$ of degree $a$.
By the work of Moeglin and Waldspurger (\cite{MW}), $\Sigma$ gives rise to a
selfdual representation in the residual spectrum of $\GL_n({\A}_K)$ denoted by $\Sigma[b]$.
Next, by the global Jacquet-Langlands correspondence, the representation
$\Sigma[b]$ of $\GL_n({\mathbb A}_K)$ can be transported  to an
automorphic representation  $\Sigma'[b]$ of $\D^\times({\mathbb A}_K)$ (of multiplicity one),
where $\D$ is a division algebra of index $n$ over $K$, which is ramified only at $u, v$, with $\D_v=D$.
Again, by Theorem 2.1, the corresponding representation of $\Sigma'[b]$ is
orthogonal at $v$ if and only if it is so at $u$. Since by construction,
$\pi'\simeq \Sigma'[b]_v$, we get $c(\pi')=c(\Sigma'[b]_u)$. Put $B:=\D_u$ and $\pi_1=
\Sigma'[b]_u$. The proof is completed by noting that an
irreducible representation $\pi$ of $D^\times/D^\times(1)$ of dimension $f$ dividing  $ n$
with $e= n/f$ has
Langlands parameter which is $ {\rm Ind}_{k_f^\times}^{W_k} (\chi) \otimes sp_e$ where
$\chi$ is a certain character of $k_f^\times$, and that Theorem A
is true for such representations of $D^\times$.

The proof of Theorem A is now complete for $D^\times$.

\qed

\bigskip

\section{\bf Proof of Theorem A for selfdual representations on $\GL_m(D)$}

\medskip

In this section we extend the results above for $D^\times$ to the case of selfdual irreducible
discrete series representations $\pi'$ of $\GL_m(D)$, with $m>1$. Put $n=md$, where $d$ is the index of $D$.
We need to prove the following:

\begin{thm} An irreducible, selfdual, discrete series
representation $\pi'$ of $\GL_m(D)$, where $D$ is a
division algebra over a non-archimedean local field $k$ of index $d$, is orthogonal if $d$ is odd.
If $d$ is even, and $m$ is odd, then the representation is orthogonal if and
only if its parameter is symplectic. If both $m$ and $d$ are even, then the
representation is orthogonal.
\end{thm}

{\it Proof}. \,
Assume first that the corresponding representation $\pi$ of GL$_n(k)$ is supercuspidal.

Suppose $d$ is odd. As $\pi$ is selfdual and supercuspidal,
we may globalize it (by applying Theorem 1.5) to
a selfdual, cuspidal automorphic representation $\Pi$ of $\GL_{n}({\mathbb A}_K)$ with $\pi$
as its local component at $K_{v} = k$. Let $E$ be a cyclic extension of
$K$ of degree $d$ such that
$v$ splits into $d$ places. Let $\Pi_E$ denote the base change of $\Pi$ to $E$.
Let $\B$ be a central division algebra over $E$ of
dimension $d^2$ over $E$ such that  $\B \otimes_K K_v \cong D^n$, and such that
$\B$ has no other ramification. That there is such a division algebra $\B$ follows
from class field theory.

Let $\Pi_E^\B$ denote the automorphic representation of $\GL_m(\B(\A_E))$ obtained
from $\Pi_E$ by the Jacquet-Langlands correspondence (\cite{bad}). The component of
$\Pi_E^D$ at $E_v:=E\otimes_K K_v \simeq k^d$ is isomorphic to ${\pi'}^{\otimes d}$.
At every place $u$ of $E $ not lying over $v$, the $u$-component
of $\Pi_E^\B$ is a representation of $\GL_{n}(E_u)$.
Applying Theorem 2.1, we find that the $d$-th power of $c(\pi')$ is trivial, and therefore
$\pi'$ is an orthogonal representation (as $d$ is odd).

Now let the index of $D$ be $d=2r$. In this case
there is a division algebra $\D$ over a number field $K$ of index $2mr$ such that
$\D$ gives rise to $D$ at one place $v$, say, necessarily with index, say
$\frac{m}{2mr} \in {\mathbb Q}/{\mathbb Z}$ there,
and with indices $\frac{-1}{2mr} \in {\mathbb Q}/{\mathbb Z}$
at $m$ other places, call them $u_1, \dots, u_m$. We may further take $\D$ to be split
at all the remaining places. The existence of such a
 division algebra follows from class field theory. This global division algebra
becomes split, or totally ramified, at every place other than $v$. Since one understands
the parity at the split and totally ramified places, the parity at the unique remaining
place $v$ follows, thanks to the global sign being $1$ by Theorem 2.1. We are now
done when $\pi$ is supercuspidal.

Now assume that $\pi$ is not supercuspidal, with its parameter being of the form
$\sigma=\tau\otimes sp_b$, $n=ab$, with $\tau$ an irreducible, selfdual
$a$-dimensional representation of $W_k$.
The proof now proceeds as in the last paragraph of section 5 except for the additional
subtlety that the global Jacquet-Langlands correspondence might produce
a Speh-like representation on $\GL_n(D)$ whereas we want to make conclusions
about discrete series representations of $\GL_n(D)$. To deal with this possibility, we note the
following proposition,
which, when combined with the fact that the Aubert-Zelevinsky
involution
interchanges Speh-modules with the generalized Steinberg
representations,
completes the proof of Theorem A. This involution $\pi\mapsto {\it i}(\pi)$
(cf. \cite{bad}, section 2.6, for example)
is defined on the Grothendieck group of
smooth representations of a $p$-adic reductive group $G$
as an alternating sum of parabolically induced
representations of the various Jacquet modules of $\pi$, and sends an irreducible
to another irreducible $\vert {\it i}(\pi)\vert$ up to sign.

\begin{prop}
Let $G$ be a reductive algebraic group over a non-archimedean
local field $k$.
Then an irreducible representation $\pi$ is orthogonal if and only if
the corresponding irreducible representation $\vert {\it i}(\pi)\vert$, associated by the
Aubert-Zelevinsky involution ${\it i}$, is orthogonal.
\end{prop}

The proof of this proposition follows from the following lemma using the
fact that both induction and Jacquet functor take real representations
to real representations.

\begin{lemma} An irreducible unitary representation $(\pi,V)$ of
a $p$-adic group $G$ carries a nonzero symmetric bilinear form
$B: V \times V \rightarrow {\Bbb C}$ if and only if $\pi$ is defined
over ${\Bbb R}$, i.e., there is a $G$-invariant real subspace $W$
of $V$ such that $V = W \otimes_{\Bbb R}{\Bbb C}$.
\end{lemma}

This lemma is well-known for finite groups, and most of the
proofs in the literature work  in this infinite
dimensional context too (assuming only that Schur's lemma holds, which
of course it does for $p$-adic groups).

\section{\bf Questions on $\theta$-selfdual representations}

\medskip

Let $G$ be a group, and $\theta$ an automorphism of $G$ with $\theta^2 =1$.
Let $\pi$ be an irreducible representation of $G$. In this section
the group $G$ will be any one of the finite, real, $p$-adic, or adelic group,
and the representations, as well as their contragredients,
will be understood in the corresponding
category. Recall that $\pi$ is called $\theta$-selfdual if $\pi^\vee
\cong \pi^{\theta}$. In section 2, we have associated an invariant
$c_\theta(\pi) \in \{ \pm 1 \}$ for any representation $\pi$ of $G$ with
$\pi^\vee \cong \pi^{\theta}$, which for $\theta = 1$ furnishes
the information on whether $\pi$ is orthogonal or symplectic.

\vspace{2mm}

\noindent{\bf Example :} Let $E$ be a quadratic extension of a local field $F$
with $\sigma$ as the Galois automorphism of $E$ over $F$. For any algebraic group
$G$ over $F$, this gives an involution, call it $\theta$, on $G(E)$, and thus
we have associated an
invariant $c_\theta\in \{ \pm 1 \}$ for any representation $\pi$ of $G(E)$ with
$\pi^\vee \cong \pi^{\theta}$.

\vspace{2mm}

\noindent{\bf Remark}: \, If $G$ is an algebraic group over a finite field ${\mathbb F}$,
${\mathbb E}$ a quadratic extension of ${\mathbb F}$,
and $\theta$ the automorphism of order $2$
on $G({\mathbb E})$ induced by the Galois action on ${\mathbb E}$, then the invariant
$c_\theta$ has been studied in \cite{P2}. It follows from the results there that for
irreducible representations $\pi$ of $G({\mathbb E})$ with
$\pi^\vee \cong \pi^{\theta}$, the invariant $c_\theta(\pi) =1$ for those
irreducible representations
$\pi$ which are linear combinations of Deligne-Lusztig representations.

\medskip

\begin{question} {\bf (a):} Let $G$ be a reductive p-adic group, and $\pi_1$
and $\pi_2$  two irreducible tempered representations of $G$. Is it
true that if $\pi_1^\vee \cong \pi_1^\theta, \pi_2^\vee \cong
\pi_2^\theta$, and $\pi_1$ and $\pi_2$ are in the same $L$-packet,
then $c_\theta(\pi_1) = c_\theta(\pi_2)$?

{\bf (b): }
If part (a) is true,
is there an expression for $c_\theta(\pi)$ in terms of the Langlands parameter of
$\pi$? The results in this
paper suggest that $c_\theta(\pi)$ depends not only on the $L$-group, but also the
inner form (as the result is different for $\GL(n)$ and for division algebras), and
that it varies as one varies the inner forms as in Kottwitz's  paper \cite{kottwitz}.

{\bf (c): } If $\pi$ is automorphic with $\pi^\vee \cong \pi^\theta$, is $c_\theta(\pi) = 1$
(like it is for $\pi$ occurring discretely with multiplicity one; cf. Theorem 2.1)?
\end{question}

When $G$ is SL$(n)$, the
elements of any $L$-packet of representations
are permuted transitively by conjugation by elements of GL$(n)$, and therefore,
at least in this case, whether a representation of SL$(n)$ is orthogonal
or symplectic is a property of the $L$-packet.

\medskip

Here is a local result concerning generic, square-integrable, $\theta$-selfdual representations in the style of Theorem 2.1.

\begin{prop} Let $G$ be a quasi-split reductive group over a local field, with maximal unipotent subgroup $N$, equipped
with an involution $\theta$ which preserves $N$. Let $\psi$ be a non-degenerate character of $N$
satisfying $\psi^\theta=\overline \psi$. Then for any $\psi$-generic discrete series representation $\pi$ of $G$ which
is $\theta$-selfdual, we have $c_\theta(\pi) = 1$.
\end{prop}

{\it Proof}, \, For any tempered $\psi$-generic $\pi$, one has an infinitesimal embedding $\xi\mapsto W_\xi$, of $\pi$ into
$L^2(N\backslash G, \psi)$, and the square-integrability of $\pi$ implies, by a theorem of Harish Chandra, that $\pi$
occurs discretely; denote by ${\mathcal W}_\pi$ the image of $\pi$. Moreover, one knows that ${\mathcal W}_\pi$ has
multiplicity one. As $\psi^\theta=\overline \psi$ by hypothesis, we may define a bilinear form $B$ on $\pi$ by
$$
(\xi, \xi') \, \rightarrow \, \int_{N\backslash G} W_\xi(g)W_{\xi'}^\theta(g) dg.
$$
Since $\pi^\theta\simeq \pi^\vee$, the function $\overline W_\xi^\theta$ lies, by multiplicity one, in ${\mathcal W}_\pi$.
So we may take $W_\xi'$ to be $\overline W_\xi^\theta$ and find that $B(\xi, \xi')$ is the $L^2$-norm of $W_\xi$. The
argument is now completed as in the proof of Theorem 2.1.

\medskip

\noindent{\bf 7.3 \, Dual involutions and dual signs}

Now suppose that $G$ is a connected reductive group over a local field $k$, with the Langlands dual group ${}^LG$ admitting a
an involution ${}^L\theta$. Such an involution on ${}^LG$ certainly exists when $G$ is quasi-split admitting an involution
$\theta$ preserving a Borel subgroup $B$, a maximal torus $T\subset B$, and the
pinned root system $\Psi$.

Consider any $L$-parameter
$$
\sigma: W_k \, \rightarrow \, {}^LG,
$$
not lying in a Levi subgroup of ${}^LG$. Suppose $\sigma$ is conjugate to $\sigma^{{}^L\theta}$ by an element $t$ of
the connected component $G^\vee$ of ${}^LG$. Then we can associate a ``sign'' $c_{_{{}^L\theta}}(\sigma)$
as follows: It is easy to see that $tt^\theta$ must belong to the center ${}^LZ$ of ${}^LG$.
This implies in particular that $t$ and $t^\theta$ commute, and that
$tt^\theta$ is a $\theta$-invariant element of ${}^LZ$. The element $t$ is
well defined modulo $^LZ$, and hence $tt^\theta$ is well-defined in the
Tate cohomolgy group ${\hat H}^0(C, {}^LZ) = {}^LZ^C/N$, where $C$ is the cyclic group $\{1, \theta\}$, ${}^LZ^C$ is
the group of $C$-fixed points of ${}^LZ$, and
$N$ the subgroup of norms. We will write $c_{_{{}^L\theta}}(\sigma)$ for the class of $tt^\theta$. If $\sigma$ is
the parameter of $\pi$, it will be interesting to compare $c_\theta(\pi)$ with $c_{_{{}^L\theta}}(\sigma)$.
When $G$ is GL$(n)/k$ with $\theta(g)={}^tg^{-1}$, we have $Z^C = \{\pm 1\}$ and $N=\{1\}$, so
$c_{_{{}^L\theta}}(\sigma)$ takes values in $\{\pm 1\}$.

\bigskip

\section{\bf Rationality}

\medskip

The question about selfdual representations being orthogonal or
symplectic is part of the more general question about {\it field of
definition} of a representation. For example, by (the well known) lemma 6.3, a selfdual,
unitary representation is orthogonal if and only if
it is defined over ${\Bbb R}$.

Let $G$ be a group, and $\pi$ an irreducible representation
of $G$ over ${\mathbb C}$. Put
$$ {\mathcal G}_\pi = \{\sigma \in {\rm Aut}({\mathbb C}/{\mathbb Q})|
\pi^\sigma \cong \pi\}.$$
If either $G$ is finite, or $G$ is a reductive $p$-adic group, and $\pi$
is supercuspidal with finite order central character,
then ${\mathcal G}_\pi$ is known to be a subgroup of
finite index of Aut$({\mathbb C}/{\mathbb Q})$, and this defines a
finite extension $K$ of ${\mathbb Q}$. Call $K$ the field
of definition  of $\pi$. (If $\pi$ is finite dimensional, it is the field
generated by the character values of $\pi$.)

Assume that $\pi$ is finite dimensional.
Then associated to $\pi$, there is a division
algebra ${\mathcal D}_\pi$ with center $K$ which
measures the obstruction to $\pi$ being defined over $K$,  called Schur Algebra, thus
${\mathcal D}_\pi = K$ if and only if $\pi$ can be defined
over $K$. Let the dimension of ${\mathcal D}_\pi$ over $K$ be $d^2$;
the integer $d$ is called the Schur index of $\pi$.

Now let $\sigma \rightarrow \pi_\sigma$ be the local Langlands
correspondence between irreducible representations of
Gal$(\bar{k}/k)$ of dimension $n$, and representations of
$\GL_m(D)$ where $D$ is a division algebra of index $r$ with $rm=n$.
Here and in what follows, we normalize the Langlands correspondence
by multiplying by the character $x\rightarrow |x|^{(n-1)/2}$ where
$x \in k^\times$. This normalized Langlands correspondence, which
does not affect rationality (reality!) questions at infinity, is
what is Galois equivariant on the coefficients; see for example, Henniart
\cite{He2}.
\vspace{3mm}

\noindent{\bf Question:} How are
${\mathcal D}_\sigma$ and ${\mathcal D}_{\pi_\sigma}$
related? This time the answer
a priori might depend not just on the index of the division algebra,
but on its class in the Brauer group. However, we propose the following
conjecture, suggesting in particular that this is not the case.

\vspace{3mm}

\begin{conj}Let $D$ be a division algebra of index $r$ over
a non-archimedean local field $k$, $\pi$ an irreducible representation
of $\GL_m(D)$, and $\sigma_{\pi}$ the associated $n$-dimensional
representation of the Weil-Deligne group of $k$ for $n=mr$.
Let ${\mathcal D}_\sigma$ and ${\mathcal D}_{\pi_\sigma}$ be the
associated Schur algebras with center a number field $K$ (which is a
cyclotomic field). Then the
following happens:

\begin{enumerate}

\item If $\pi$ is not selfdual, or $\pi$ is selfdual with
$c(\sigma) = c(\pi_\sigma)$, then
${\mathcal D}_\sigma = {\mathcal D}_{\pi_\sigma}$.

\item If $\pi$ is selfdual and $c(\pi) = -c(\pi_\sigma)$,
then the answer depends
on the degree  $K$ over ${\mathbb Q}$ which is a certain
totally real extension of ${\mathbb Q}$.
If $[K:{\mathbb Q}]$
is even, then the invariants of ${\mathcal D}_\sigma$
and ${\mathcal D}_{\pi_\sigma}$ are the same except at infinite places,
where the invariants of ${\mathcal D}_\sigma$ and ${\mathcal D}_{\pi_\sigma}$
differ by $1/2$. If $[K:{\mathbb Q}]$
is odd, then in particular there are odd number of places in $K$
over $p$. The invariants of ${\mathcal D}_\sigma$
and ${\mathcal D}_{\pi_\sigma}$ are the same except at infinite places, and
at the places in $K$ above $p$,
where the invariants of ${\mathcal D}_\sigma$ and ${\mathcal D}_{\pi_\sigma}$
differ by $1/2$.
\end{enumerate}

\end{conj}

One case of the conjecture is especially simple to state. This is
when $r$, the index of $D$  is odd, and $\pi$ is a selfdual representation
of $\GL_m(D)$ for $m$ odd. In this case,
irreducible selfdual representations of the
Galois group of dimension $n$, or irreducible selfdual
representation of $\GL_m(D)$ exist only in even residue characteristic,
cf. \cite{P1}. We are thus in the tame case $(n,p)=1$, and in this case
it can be seen that  Galois representations are induced from a
character $\theta$ of $L^\times$ where $L$ is a degree $n$ extension of
$k$, with $\theta^2 =1$, thus $\theta$ takes values in $\pm 1$.
Therefore in this case the Galois representation is defined over
${\mathbb Q}$. Our theorem B implies  that the selfdual
representations of $\GL_m(D)$ are defined over ${\mathbb R}$, and the
discussion in this section refines it to ask the following:

\medskip

\noindent{{\bf Question 8.2} \, {\it Is every selfdual, irreducible
representation of $\GL_m(D)$, for $D/k$ a division algebra
of odd index, defined over ${\mathbb Q}$ for $m$ odd?}

\medskip

Recently, Bushnell and Henniart have answered this question in the affirmative
in \cite{BH} for $m=1$ or $r =1$.

\bibliographystyle{math}    
\bibliography{selfdualreps}

\def\cftil#1{\ifmmode\setbox7\hbox{$\accent"5E#1$}\else
  \setbox7\hbox{\accent"5E#1}\penalty 10000\relax\fi\raise 1\ht7
  \hbox{\lower1.15ex\hbox to 1\wd7{\hss\accent"7E\hss}}\penalty 10000
  \hskip-1\wd7\penalty 10000\box7}
\begin{thebibliography}{CKPSS}

\bibitem[AC]{AC}
J.~Arthur and L.~Clozel.
\newblock {\em Simple algebras, base change, and the advanced theory of the
  trace formula}, volume 120 of {\em Annals of Mathematics Studies}.
\newblock Princeton University Press, Princeton, NJ, 1989.

\bibitem[Bad]{bad}
A.~I. Badulescu.
\newblock {Global {J}acquet-{L}anglands correspondence, multiplicity one and
  classification of automorphic representations}.
\newblock {\em Invent. Math.} {\bf 172} (2008), 383--438.
\newblock With an appendix by Neven Grbac.

\bibitem[BH]{BH}
C.~Bushnell and G.~Henniart.
\newblock {Self-dual representations of some dyadic groups}, 2008.

\bibitem[CC]{Ch-Cl}
G.~Chenevier and L.~Clozel.
\newblock {Corps de nombres peu ramifies et formes automorphes autoduales},
  2007.

\bibitem[Clo1]{Cl2}
L.~Clozel.
\newblock {On limit multiplicities of discrete series representations in spaces
  of automorphic forms}.
\newblock {\em Invent. Math.} {\bf 83} (1986), 265--284.

\bibitem[Clo2]{Cl}
L.~Clozel.
\newblock {Repr\'esentations galoisiennes associ\'ees aux repr\'esentations
  automorphes autoduales de {${\rm GL}(n)$}}.
\newblock {\em Inst. Hautes \'Etudes Sci. Publ. Math.} (1991), 97--145.

\bibitem[CKPSS]{CoKPSS}
J.~W. Cogdell, H.~H. Kim, I.~I. Piatetski-Shapiro, and F.~Shahidi.
\newblock {Functoriality for the classical groups}.
\newblock {\em Publ. Math. Inst. Hautes \'Etudes Sci.} (2004), 163--233.

\bibitem[CPSS]{CoPSS}
J.~W. Cogdell, I.~I. Piatetski-Shapiro, and F.~Shahidi.
\newblock {Stability of {$\gamma$}-factors for quasi-split groups}.
\newblock {\em J. Inst. Math. Jussieu} {\bf 7} (2008), 27--66.

\bibitem[GJR]{GJR}
D.~Ginzburg, D.~Jiang, and S.~Rallis.
\newblock {On the nonvanishing of the central value of the {R}ankin-{S}elberg
  {$L$}-functions}.
\newblock {\em J. Amer. Math. Soc.} {\bf 17} (2004), 679--722 (electronic).

\bibitem[HL]{HaLab}
M.~Harris and J.-P. Labesse.
\newblock {Conditional base change for unitary groups}.
\newblock {\em Asian J. Math.} {\bf 8} (2004), 653--683.

\bibitem[HT]{HaT}
M.~Harris and R.~Taylor.
\newblock {\em The geometry and cohomology of some simple {S}himura varieties},
  volume 151 of {\em Annals of Mathematics Studies}.
\newblock Princeton University Press, Princeton, NJ, 2001.
\newblock With an appendix by Vladimir G. Berkovich.

\bibitem[Hen1]{He}
G.~Henniart.
\newblock {Une preuve simple des conjectures de {L}anglands pour {${\rm
  GL}(n)$} sur un corps {$p$}-adique}.
\newblock {\em Invent. Math.} {\bf 139} (2000), 439--455.

\bibitem[Hen2]{He2}
G.~Henniart.
\newblock {Une caract\'erisation de la correspondance de {L}anglands locale
  pour {${\rm GL}(n)$}}.
\newblock {\em Bull. Soc. Math. France} {\bf 130} (2002), 587--602.

\bibitem[HH]{HH}
G.~Henniart and R.~Herb.
\newblock {Automorphic induction for {${\rm GL}(n)$} (over local
  non-{A}rchimedean fields)}.
\newblock {\em Duke Math. J.} {\bf 78} (1995), 131--192.

\bibitem[JS1]{J-S1}
D.~Jiang and D.~Soudry.
\newblock {The local converse theorem for {${\rm SO}(2n+1)$} and applications}.
\newblock {\em Ann. of Math. (2)} {\bf 157} (2003), 743--806.

\bibitem[JS2]{J-S2}
D.~Jiang and D.~Soudry.
\newblock {Generic representations and local {L}anglands reciprocity law for
  {$p$}-adic {${\rm SO}\sb {2n+1}$}}.
\newblock In {\em Contributions to automorphic forms, geometry, and number
  theory}, pages 457--519. Johns Hopkins Univ. Press, Baltimore, MD, 2004.

\bibitem[Kot]{kottwitz}
R.~E. Kottwitz.
\newblock {Sign changes in harmonic analysis on reductive groups}.
\newblock {\em Trans. Amer. Math. Soc.} {\bf 278} (1983), 289--297.

\bibitem[Lab]{Lab}
J.-P. Labesse.
\newblock {Cohomologie, stabilisation et changement de base}.
\newblock {\em Ast\'erisque} (1999), vi+161.
\newblock Appendix A by Laurent Clozel and Labesse, and Appendix B by Lawrence
  Breen.

\bibitem[MW]{MW}
C.~M{\oe}glin and J.-L. Waldspurger.
\newblock {Le spectre r\'esiduel de {${\rm GL}(n)$}}.
\newblock {\em Ann. Sci. \'Ecole Norm. Sup. (4)} {\bf 22} (1989), 605--674.

\bibitem[M{\oe}g]{Moeglin}
C.~M{\oe}glin.
\newblock {Classification et changement de base pour les s\'eries discr\`etes
  des groupes unitaires {$p$}-adiques}.
\newblock {\em Pacific J. Math.} {\bf 233} (2007), 159--204.

\bibitem[Pra1]{P2}
D.~Prasad.
\newblock {Distinguished representations for quadratic extensions}.
\newblock {\em Compositio Math.} {\bf 119} (1999), 335--345.

\bibitem[Pra2]{P1}
D.~Prasad.
\newblock {Some remarks on representations of a division algebra and of the
  {G}alois group of a local field}.
\newblock {\em J. Number Theory} {\bf 74} (1999), 73--97.

\bibitem[PR]{PR}
D.~Prasad and D.~Ramakrishnan.
\newblock {Lifting orthogonal representations to spin groups and local root
  numbers}.
\newblock {\em Proc. Indian Acad. Sci. Math. Sci.} {\bf 105} (1995), 259--267.

\bibitem[Rog]{Rog}
J.~D. Rogawski.
\newblock {\em Automorphic representations of unitary groups in three
  variables}, volume 123 of {\em Annals of Mathematics Studies}.
\newblock Princeton University Press, Princeton, NJ, 1990.

\bibitem[SZ]{S-Z}
A.~J. Silberger and E.-W. Zink.
\newblock {An explicit matching theorem for level zero discrete series of unit
  groups of {$\mathfrak p$}-adic simple algebras}.
\newblock {\em J. Reine Angew. Math.} {\bf 585} (2005), 173--235.

\bibitem[Sou1]{Sou}
D.~Soudry.
\newblock {On {L}anglands functoriality from classical groups to {${\rm GL}\sb
  n$}}.
\newblock {\em Ast\'erisque} (2005), 335--390.
\newblock Automorphic forms. I.

\bibitem[Sou2]{Sou-lett}
D.~Soudry.
\newblock {Letter to D.~Ramakrishnan}, December 2007.

\bibitem[Vig]{V}
M.-F. Vign{\'e}ras.
\newblock {Correspondances entre representations automorphes de {${\rm GL}(2)$}
  sur une extension quadratique de {${\rm GSp}(4)$} sur {${\bf Q}$}, conjecture
  locale de {L}anglands pour {${\rm GSp}(4)$}}.
\newblock In {\em The Selberg trace formula and related topics (Brunswick,
  Maine, 1984)}, volume~53 of {\em Contemp. Math.}, pages 463--527. Amer. Math.
  Soc., Providence, RI, 1986.

\end{thebibliography}

\vspace{0.6in}

\noindent
\begin{tabular}{l}
School of Mathematics,
Tata Institute of Fundamental Research,
Mumbai- 400 005, India.\\
Email: dprasad@math.tifr.res.in\\

\\

Department of Mathematics,
California Institute of Technology,
Pasadena, CA 91125.\\
Email: dinakar@caltech.edu

\end{tabular}

\end{document}